\documentclass[12pt,a4paper]{article}
\usepackage{amssymb,amsmath,amsthm,color}
\usepackage{url}
\newcommand{\bt}{\begin{Theorem}}
\newcommand{\et}{\end{Theorem}}
\newcommand{\bl}{\begin{Lemma}}
\newcommand{\el}{\end{Lemma}}
\newcommand{\br}{\begin{Remark}}
\newcommand{\er}{\end{Remark}}

\newcommand{\R}{\mbox{$\mathbb R$}}
\newcommand{\C}{\mbox{$\mathbb C$}}

\newcommand{\bi}{\begin{itemize}}
\newcommand{\ei}{\end{itemize}}
\newcommand{\bea}{\begin{eqnarray}}
\newcommand{\eea}{\end{eqnarray}}

\newtheorem{Definition}{Definition}[section]
\newtheorem{Theorem}[Definition]{Theorem}
\newtheorem{Lemma}[Definition]{Lemma}

\newtheorem{Remark}[Definition]{Remark}

\newcommand{\be}{\begin{equation}}
\newcommand{\ee}{\end{equation}}

\theoremstyle{plain} \theoremstyle{plain}
\newtheorem{theorem}{Theorem}[section]
\newtheorem{lemma}[theorem]{Lemma}

\theoremstyle{definition}

\theoremstyle{plain}
\begin{document}
\begin{center}
{\large \bf { Twisted spherical  means in annular regions \\in
$\mathbb C ^n$ and support theorems}}
\end{center}
\begin{center}
{\bf Rama Rawat  and  R. K. Srivastava}
\end{center}
\begin{center}
{\bf Dedicated to Prof. A. Sitaram on his sixtieth birthday}
\end{center}

\begin{abstract}
Let $Z(Ann(r,R))$ be the class of all continuous functions $f$ on
the annulus $Ann(r,R)$ in $\mathbb C^n$ with twisted spherical
mean $f \times \mu_s(z)=0,$
whenever $z\in \mathbb C^n$ and $s >0$ satisfy the condition that the sphere
$S_s(z)\subseteq Ann(r, R) $ and ball $B_r(0)\subseteq B_s(z).$
In this paper, we give a
characterization for functions in $Z(Ann(r,R))$ in terms of their
spherical harmonic coefficients. We also prove support theorems for the twisted
spherical means in $\mathbb C^n$ which improve some of the earlier
results.
\end{abstract}
\begin{center}
\section{\bf Introduction and the main results}
\end{center}
 For $s >0,$ let $\mu_s$ stand for the normalized
surface measure on $\{z \in \mathbb C^n : |z| = s \}.$ The twisted
spherical means of a function $f$ in $L^1_{loc}(\mathbb C^n)$ are
defined by

\begin{equation}\label{exp1}
f \times \mu_s (z) = \int_{|w| = s }~
f(z-w)e^{\frac{i}{2}Im(z.\bar{w})} d\mu_s(w), ~~ z\in \mathbb C^n.
\end{equation}
 These twisted spherical means arise in a natural way from the
 spherical means on the Heisenberg group $\mathbb H^n .$  The
 group $\mathbb H^n,$ as a manifold is $\mathbb C^n \times
\mathbb R$, with the group law $$(z, t) (w, s) = (z+w,
t+s+\frac{1}{2}Im z.\bar{w}).$$

If $\mu_s$ is now considered as a measure on $ \{ (z, 0):~|z| = s
\} \subset \mathbb H^n,$ then the spherical means of a function
$f$ in $L^1_{loc}(\mathbb H^n)$ are defined by
\begin{equation}\label{exp2}
f \ast \mu_s (z, t) = \int_{|w|=s}~f((z,t)(-w,0))~d\mu_s(w).
\end{equation}

Let
$$f^\lambda(z)= \int_\mathbb R f(z,t)e^{i \lambda t} dt,$$
 be the inverse Fourier transform of $f$ in the $\mathbb R$
variable.

Then a simple calculation shows that
\begin{equation}\label{exp3}
(f \ast \mu_s)^\lambda=\int_{-\infty}^{~\infty}~f \ast \mu_s(z,t)
e^{i\lambda t} dt = \int_{|w| = s}~f^\lambda (z-w)
e^{\frac{i\lambda}{2} Im(z.\bar{w})}~d\mu_s(w).
\end{equation}

We can also define the $\lambda$-twisted convolution of functions
$F$ and $G$ in $L^1(\mathbb C^n)$  by
$$ F \times_\lambda G(z) = \int_{\mathbb
C^n}~F(z-w)~G(w)~e^{\frac{i\lambda}{ 2}Im(z.\bar{w})}~dw.$$

Then, (\ref{exp3}) can be rewritten as $$ (f \ast \mu_s)^\lambda
(z) = f^\lambda \times_\lambda \mu_s(z).$$

Thus, the spherical means $f \ast \mu_s$ on the Heisenberg group
can be studied using the $\lambda$-twisted spherical means
$f^\lambda \times_\lambda \mu_s$ on $\mathbb C^n.$ A further
scaling argument shows that it is enough to study these means for
the case of $\lambda = 1.$ From now onwards, we shall write $F
\times G$ instead of $F\times_1 G$ and call it the twisted
convolution of $F$ and $G.$

Let $Ann(r, R)=\{z\in \mathbb C^n :~r<|z|<R\},$ $ 0\leq r<R\leq
\infty,$ be an open annulus in $\mathbb C^n.$ Let $Z(Ann(r, R))$
be the class of all continuous functions on $Ann(r, R)$ with the
twisted spherical means
$$
 \int_{|w| = s }~
f(z-w)e^{\frac{i}{2}Im(z.\bar{w})} d\mu_s(w)=0,$$ for all $z\in \mathbb C^n$  and $s>0$ satisfying the condition that the sphere
$S_s(z)$ is contained in the annulus $Ann(r, R)$ and the ball $B_s(z)$ contains the ball $B_r(0).$

Equivalently, $f \in Z(Ann(r, R))$ if  $f \times \mu_s(z)=0,$ for all $z \in \mathbb C^n$ and $s>0$ for which the sphere
$S_s(z)$ is contained in the annulus $Ann(r, R)$ and the ball $B_s(z)$ contains the ball $B_r(0).$

Our main result, Theorem $1.1$, gives a necessary and sufficient
condition for a function $f$ to be in $Z(Ann(r,R))$ in terms of
its spherical harmonic coefficients. As a corollary, we shall also
prove some support theorems, for the twisted spherical means,
which improve results in \cite{NT2}.

This work is motivated, in spirit, by the work of Epstein and
Kliener \cite{EK} on the spherical means in annular regions in
Euclidean spaces. For some other closely related work on spherical means see \cite{
AR}, \cite{NT1}.

To state our results, we shall require the following basic facts
from the theory of bigraded spherical harmonics. (See \cite{T}, p. 12).
We shall use the notation $K=U(n)$ and $M=U(n-1).$ Then
$S^{2n-1}\cong K/M$ under the map $kM\rightarrow k.e_n,$ $k\in U(n)$
where $e_n=(0,0,\ldots ,1)\in \mathbb C^n.$ Let $\hat{K}_M$ denote
the set of all the equivalence classes of irreducible unitary
representations of $K$ which have a nonzero $M$-fixed vector. For
our set up of $K$ and $M$, it is known that for each representation
in $\hat{K}_M$ has a unique nonzero $M$-fixed vector, up to a scalar
multiple.

For a $\delta \in \hat{K}_M,$ which is realized on $V_{\delta},$
let $\{e_1,\ldots ,e_{d(\delta)}\}$ be an orthonormal basis of
$V_{\delta},$ with $e_1$ as the $M-$ fixed vector. Let
$t_{ij}^{\delta}(k)=\langle e_i,\delta (k)e_j \rangle ,$ $k\in K$
and $\langle,\rangle $ stand for the innerproduct on $V_{\delta}.$
By Peter-Weyl theorem, it follows that $\{\sqrt{d(\delta
)}t_{j1}^{\delta}:1\leq j\leq d(\delta ),\delta\in\hat{K}_M\}$ is
an orthonormal basis of $L^2(K/M).$ (see \cite{T}, p. 14 for
details). Define $Y_j^{\delta} (\omega )=\sqrt{d(\delta
)}t_{j1}^{\delta}(k),$ where $\omega =k.e_n\in S^{2n-1},$ $k \in
K.$ It then follows that $\{Y_j^{\delta}:1\leq j\leq d(\delta
),\delta\in \hat{K}_M, \}$ forms an orthonormal basis for
$L^2(S^{2n-1}).$

For our purposes, we need a concrete realization of the
representations in $\hat{K}_M,$ which can be done in the following
way. See \cite{R}, p. 253, for details.

Let $\mathbb Z^+$ denote the set of all non negative integers. For
$p, q \in \mathbb Z^+$, let $P_{p,q}$ denote the space of all
polynomials $P$ in $z$ and $\bar{z}$ of the form
$$ P(z) = \sum_{|\alpha| = p} \sum_{|\beta| = q}c_{\alpha
\beta}~z^\alpha \bar{z}^\beta .$$ Let $H_{p,q} = \{ P \in
P_{p,q}~:~\Delta P =0 \}$ where $\Delta$ is the standard Laplacian
on $\mathbb C^n.$ The elements of $H_{p,q}$ are called the bigraded
solid harmonics on $\mathbb C^n.$ The group  $K$ acts on $H_{p,q}$
in a natural way. It is easy to see that the space $H_{p,q}$ is
$K$-invariant. Let $\pi_{p,q}$ denote the corresponding
representation of $K$ on $H_{p,q}.$ Then, representations in
$\hat{K}_M$ can be identified, up to unitary equivalence, with the
collection $\{\pi_{p,q}: p,q \in \mathbb Z^+.\}$

Define the bigraded spherical harmonics on the sphere $S^{2n-1}$
by $Y_j^{p,q}(\omega )=\sqrt{d(p,q )}t_{j1}^{p,q}(\sigma ),$ where
$\omega =k.e_n\in S^{2n-1},$ $k \in K$ and $d(p,q)$ is the
dimension of $H_{p,q}.$  Then $\{Y_j^{p,q}:1\leq j\leq d(p,q),p,q
\in \mathbb Z^+ \}$ forms an orthonormal basis for
$L^2(S^{2n-1}).$

Therefore, for a continuous function $f$ on $\mathbb C^n,$ writing
$z=\rho \,\omega,$ where $\rho
> 0$ and $\omega \in S^{2n-1},$ we can expand the function $f$ in
terms of spherical harmonics as
\begin{equation}\label{exp4}
f (\rho\omega) = \sum_p \sum_q \sum_{j=1}^{d(p,q)}~
a_j^{p,q}(\rho)~Y_j^{p,q}(\omega).
\end{equation}

The functions $~a_j^{p,q} $ are called the spherical harmonic
coefficients of the function $f$.

The $(p,q)^{th}$ spherical harmonic projection, $\Pi_{p,q}(f)$, of
the function $f$ is then defined as
\begin{equation}\label{exp5}
\Pi_{p,q}(f)(\rho,\omega)=~\sum_{j=1}^{d(p,q)}~a_j^{p,q}(\rho)~Y_j^{p,q}(\omega).
\end{equation}

We will replace the spherical harmonic $Y_j^{p,q}(\omega)$ on the
sphere by the solid harmonic
$P_j^{p,q}(z)=|z|^{p+q}Y_j^{p,q}(\frac{z}{|z|})$ on $\mathbb C^n$
and accordingly for a function $f$, define
$\tilde{a}_j^{p,q}(\rho)= \rho^{-(p+q)}~a_j^{p,q}(\rho),$ where
$a_j^{p,q}$ are defined by equation \ref{exp4}. We shall continue
to call the functions $\tilde{a}_j^{p,q}$ the spherical harmonic
coefficients of $f.$

Our main result is the following characterization
theorem.
\vskip.1in

\bt\label{th1}

Let $f(z)$ be a continuous function on $Ann(r, R).$ Then a
necessary and sufficient condition for $f$ to be in $Z(Ann(r,
R))$ is that for all $p,q \in \mathbb Z^+,$ $ 1\leq j \leq
d(p,q),$ the spherical harmonic coefficients $\tilde{a}_j^{p,q}$
of $f$ satisfy the following conditions:

\begin{enumerate}

\item For $p=0, q=0,$ and $r<\rho<R,$
\[~\tilde{a}_j^{0,0}(\rho)=0.\]

\item For $p,q \geq 1,$ and $r<\rho<R,$ there exists $c_i, d_k \in
\mathbb C,$ such that

$$~\tilde{a}_j^{p,q}(\rho)=\sum_{i=1}^{p}~c_i~e^{\frac{1}{4}\rho^2}\rho^{-2(p+q+n-i)}+
\sum_{k=1}^{q}~d_k~e^{-\frac{1}{4}\rho^2}\rho^{-2(p+q+n-k)}.$$

\item For $q=0$ and $p\geq 1$ or $p=0$ and $q\geq 1,$ and
$r<\rho<R,$ there exists $c_i, d_k \in \mathbb C,$ such that
$$~\tilde{a}_j^{p,0}(\rho)=~\sum_{i=1}^{p}~c_{i}~e^{\frac{1}{4}\rho^2}\rho^{-2(p+n-i)},
~\tilde{a}_j^{0,q}(\rho)=~\sum_{k=1}^{q}~d_{k}~e^{-\frac{1}{4}\rho^2}\rho^{-2(q+n-k)}.$$

\end{enumerate}
 \et
Using the above characterization for the case when $R= \infty,$ we also prove the following support
theorems for the twisted spherical means.
 \bt\label{th2}

  Let $f$ be a continuous function on $\mathbb C^n$ such that for each $k=0,
1,2,\cdots,$
$|z|^k e^{\frac{1}{4}|z|^2}|f(z)|\leq C_k.$ Then $f$ is supported in $|z| \leq r$ if and only if
$f \times \mu_s(z) = 0 $ for $s > r + |z|$ and for every $z \in
\mathbb C^n.$ \et

\bt\label{th3}

Let $f$ be a continuous function on $\mathbb C.$
 Then $f$ is supported in $|z| \leq r$
if and only if $f \times \mu_s (z) = \mu_s \times f(z)  = 0 $ for
$ s > r+ |z|$ and for every
$z \in \mathbb C.$ \et

\section{\bf Preliminaries}

We begin with the observation that the $U(n)$-invariance of the
annulus and the measure $\mu_s$ implies that for any $f$ in
$Z(Ann(r,R))$ and $p, q \in \mathbb Z^+$, $\Pi_{p,q}(f),$ as
defined in equation \ref{exp5}, also belongs to $Z(Ann(r,R)).$ In
fact the following stronger result is true.

\begin{lemma}\label{lemma1}
Suppose $f \in Z(Ann(r, R)).$ Then for $p,q \in \mathbb Z^+$,
$$a_j^{p,q}(|z|) Y_i^{p,q}(\omega) \in Z(Ann(r, R)), 1 \leq i,j
\leq d_{p,q}.$$ In particular, if ~$f \in Z(Ann(r, R)),$ then
$\Pi_{p,q}(f) \in Z(Ann(r,R))$ for all $p,q \in \mathbb Z^+.$
\end{lemma}

\begin{proof}: For $k \in U(n), \omega \in S^{2n-1},$ we have
$$Y_i^{p,q} (k^{-1} \omega ) =
\sum_{j=0}^{d(p,q)}~\overline{t_{ji}^{p,q} (k)}~
Y_j^{p,q}(\omega).$$

Using the orthogonality of the matrix entries, we have
\begin{equation}\label{exp6}
a_j^{p,q}(|z|)~Y_i^{p,q}(\omega) = d(p,q)~\int_{U(n)}~f(k^{-1} z)
t_{ij}^{p,q}(k)~dk
\end{equation} for $1 \leq i,~j \leq d(p,q).$

The proof now follows from the $U(n)$-invariance of the annulus
and the measure $\mu_s.$
\end{proof}
We shall also frequently need the following lemma to decompose a
homogeneous polynomial into sum of homogeneous harmonic
polynomials uniquely.

\begin{lemma}\label{lemma2}
Let $P\in P_{p,q}.$ Then we can write
$P(z)=P_0(z)+|z|^2P_1(z)+.........+|z|^{2l}P_l(z)$ where $P_k\in
H_{p-k, q-k}$, and $l\leq min(p,q).$
\end{lemma}

For a proof of this lemma see \cite{T}, p. 66.

Let $p, q, l, m \in \mathbb Z^+.$ Define the space $H_{p,q} \cdot
H_{l,m}$ to be the vector space of finite sums of the form $\sum
P_i Q_i$ where $P_i \in H_{p,q}$ and $Q_i \in H_{l,m}.$ Let
$$\nu = \nu(p, q, l, m) = min(p, m) + min(l,q).$$ Then the
following lemma has been proved in \cite{R}, p. 253.

\begin{lemma}\label{lemma3}
$H_{p,q}. H_{l,m} \subset \sum_{j=0}^{\nu}~H_{p+l-j,~ q+m-j} \quad
\mbox{where} \quad \nu = \nu(p, q, l,m).$
\end{lemma}

As in the proof of the Euclidean case \cite{EK}, to characterize
functions in $Z(Ann(r, R),)$ it would be enough to characterize
the spherical harmonic coefficients of smooth functions in
$Z(Ann(r, R))$. This can be done using the following approximation
argument. Let $\phi$ be a nonnegative, radial, smooth, compactly
supported function supported in the unit ball in $\mathbb C^n$
with $\int_{\mathbb C^n}\phi=1.$

Let $\phi_\epsilon(z)= \epsilon^{-2n}\phi(\frac{z}{\epsilon}).$
Then the function $$S_\epsilon(f)(z)=\int_{\mathbb C^n}
f(z-w)\phi_\epsilon(w)e^{\frac{i}{2}Im(z.\bar{w})}
dw$$ is smooth and it is easy to see that
$S_\epsilon(f)$ lies in $Z(Ann(r+\epsilon, R-\epsilon))$ for each
 $\epsilon >0.$ Since $f$ is continuous,
 $S_\epsilon(f)$ converges to $f$ uniformly on compact sets. Therefore, for each $p,q$,
 $$\lim\limits_{\epsilon \rightarrow
 0}\Pi_{p,q}(S_\epsilon(f))=\Pi_{p,q}(f).$$

Henceforth, we would assume, without loss of generality, that the functions in
$Z(Ann(r, R))$ are also smooth in the annulus $Ann(r, R).$ This
would allow us to differentiate the functions in $Z(Ann(r, R))$
arbitrarily.

Let us define the $2n$ vector fields on $\mathbb C^n$ by
$$Z_j = \frac{\partial}{\partial z_j} - \frac{1}{4} \bar{z_j},
~~~~\bar{Z_j} = \frac{\partial}{\partial \bar{z_j} } + \frac{1}{4}
z_j,~~~~j = 1, 2, \cdots \cdots n.$$

These vector fields together with the identity generate an algebra
which is isomorphic to the $(2n+1)$ dimensional Heisenberg
algebra. For the twisted convolution on $\mathbb C^n,$ they play a
role similar to that of the Lie algebra of left invariant vector
fields on a Lie group.

It is easy to verify that if $f \in Z(Ann(r, R)),$ then

$$ Z_j (f \times \mu_s) = Z_jf \times \mu_s~and~\bar{Z_j} ( f
\times  \mu_s) = \bar{Z_j}f \times  \mu_s. $$

As a consequence, $Z_j f$ and $\bar{Z_j} f$ both belong to
$Z(Ann(r, R)).$
\section{\bf The Proofs}

We shall first prove the necessary part of Theorem \ref{th1}.
For this, by Lemma \ref{lemma1}, it is enough to prove the
following theorem.

\bt\label{th4} Let $f$ be a smooth function on $Ann(r, R)$ of the
form $f(z)=\tilde{a}(\rho)~P(z),$ where $|z|=\rho$ and $ P \in
H_{p,q}.$ Then, for $f$ to be in $Z(Ann(r, R))$ it is necessary
that $\tilde{a}$ satisfies the following conditions.

\begin{enumerate}

\item If $p=0, q=0$ and $r<\rho<R,$ then $~\tilde{a}(\rho)=0.$

\item If $p,q \geq 1$ and  $r<\rho<R,$ then there exists $c_i, d_k
\in \mathbb C,$ such that
$$~\tilde{a}(\rho)=\sum_{i=1}^{p}~c_i~e^{\frac{1}{4}\rho^2}\rho^{-2(p+q+n-i)}+
\sum_{k=1}^{q}~d_k~e^{-\frac{1}{4}\rho^2}\rho^{-2(p+q+n-k)}.$$

\item If $q=0$ and $p\geq 1$ and $r<\rho<R,$ then there exists
$c_i \in \mathbb C,$ such that
$$\tilde{a}(\rho)=~\sum_{i=1}^{p}~c_{i}~e^{\frac{1}{4}\rho^2}\rho^{-2(p+n-i)}.$$

\item If $p=0$ and $q\geq 1,$ and $r<\rho<R,$ then there exists
$d_k \in \mathbb C,$ such that
$$\tilde{a}(\rho)=~\sum_{k=1}^{q}~d_{k}~e^{-\frac{1}{4}\rho^2}\rho^{-2(q+n-k)}.$$
\end{enumerate}
 \et

\begin{proof} If $p=0, q=0,$ then

\[\tilde{a}(\rho)=\int_{|w|=\rho}f(w)d\mu_\rho(w)=f\times\mu_\rho(0)=0
~\mbox{for}~ R>\rho>r,\] and the condition on $\tilde{a}_{0,0}$
follows.

For the other cases, we proceed in the following way. Since
$\bar{Z_j}f \in Z(Ann(r, R))$, computing $$\bar{Z_j}f =
\frac{\partial{f}}{\partial \bar{z_j} } + \frac{1}{4} z_jf,$$ we
have
\[\bar{Z_j}f ={\frac{z_j}{2\rho}}{\frac{\partial{\tilde{a}}}{\partial \rho }}P
 + \tilde{a}\frac{\partial{P}}{\partial \bar{z_j} } + \frac{1}{4}\tilde{a}z_jP,\]
 i.e.,
\begin{equation}
\bar{Z_j}f~=\frac{1}{2}\left(\frac{1}{\rho}\frac{\partial{\tilde{a}}}{\partial\rho}
+ \frac{1}{2}\tilde{a}\right)z_jP+
\tilde{a}\frac{\partial{P}}{\partial\bar{z_j}}~.
\end{equation}
Also
\begin{eqnarray*}
\triangle_z(z_jP)&=&4\sum_{k=1}^{n}\frac{\partial^2}{{\partial{z_k}}{\partial\bar{{z_k}}}}(z_jP)\\
                 &=&4\frac{\partial^2}{{\partial{z_j}}{\partial\bar{{z_j}}}}(z_jP)+4\sum_{k\neq
                 j}\frac{\partial^2}{{\partial{z_k}}{\partial\bar{{z_k}}}}(z_jP)\\
                 &=&4\frac{\partial{P}}{\partial \bar{z_j} } + z_j\triangle_z(P).\\
\end{eqnarray*}
Since P is harmonic, we have
\begin{equation}\label{exp7}
\triangle_z(z_jP)=~4\frac{\partial{P}}{\partial\bar{z_j}}.
\end{equation}

We shall need the identity
\begin{equation}\label{exp8}
\triangle_z(|z|^2\frac{\partial{P}}{\partial\bar{z_j}})=~4(n+p+\overline{q-1})\frac{\partial{P}}{\partial
\bar{z_j} }.
\end{equation}

For this, note that
\begin{eqnarray*}
\triangle_z\left(|z|^2P\right)&=&4\sum_{k=1}^{n}\frac{\partial^2}{{\partial{z_k}}{\partial\bar{{z_k}}}}\left(|z|^2P\right)\\
                              &=& 4\sum_{k=1}^{n}\frac{\partial}{\partial\bar{{z_k}}}\left(\bar{z_k}P
+|z|^2\frac{\partial{P}}{\partial\bar{z_k}}\right)\\
                              &=&4\sum_{k=1}^{n}\left[P+\bar{z_k}\frac{\partial{P}}{\partial\bar{z_k}}+
z_k\frac{\partial{P}}{\partial{z_k}}\right]\\
                              &=&4(n+q+p)P.
\end{eqnarray*}

Since~$\dfrac{\partial{P}}{\partial\bar{z_j}}$ is a homogeneous
harmonic polynomial of degree $ p+(q-1)$, we have (\ref{exp8}). By
Lemma \ref{lemma2},~~$z_jP(z)\in P_{p+1,q}$~ has a unique
representation
\begin{equation}\label{exp9}
z_jP(z)=P_0(z)+|z|^2P_1(z)+\dots+|z|^{2l}P_l(z)
\end{equation}
where $P_k\in H_{p+1-k,q-k}, 1 \leq k \leq l\leq \min(p+1,q).$ We
shall now show that
\begin{equation}\label{exp10}
z_jP(z)=P_0(z)+\frac{{\rho^2}}{(n+p+q-1)}\frac{\partial{P}}{\partial\bar{z_j}}.
\end{equation}
From (\ref{exp7}) and (\ref{exp8}), we have
\[\Delta_z\left(z_jP-\frac{{\rho^2}}{(n+p+q-1)}\frac{\partial{P}}{\partial\bar{z_j}}\right)=0.\]

We know that representation in (\ref{exp9}) is unique. Therefore
\[z_jP(z)=\left[z_jP(z)-\frac{{|z|^2}}{(n+p+q-1)}\frac{\partial{P}}{\partial\bar{z_j}}\right]
 +\frac{{|z|^2}}{(n+p+q-1)}\frac{\partial{P}}{\partial\bar{z_j}}\]
which is nothing but (\ref{exp10}). In view of (\ref{exp10}),
(3.7) can be rewritten as
$$\bar{Z_j}f(z)=\frac{{1}}{2}\left(\frac{{1}}{\rho}\frac{\partial{\tilde{a}}}{\partial{\rho}}+
\frac{1}{2}\tilde{a}\right)\left[P_0(z)+\frac{{\rho^2}}{(n+p+q-1)}\frac{\partial{P}}{\partial\bar{z_j}}\right]+\tilde{a}\frac{\partial{P}}{\partial
\bar{z_j}}.$$

After rearranging the terms, we have
\begin{eqnarray*}
\bar{Z_j}f(z)&=&\frac{{1}}{2}\left(\frac{{1}}{\rho}\frac{\partial{\tilde{a}}}{\partial{\rho}}+\frac{{1}}{2}\tilde{a}\right)
P_0 \\&&+
\left[\left\{\frac{{1}}{2(n+p+q-1)}\left(\rho\frac{\partial}{\partial{\rho}}+\frac{{1}}{2}\rho^2
\right)+1\right\}\tilde{a}\right]\frac{\partial{P}}{\partial\bar{z_j}}.
\end{eqnarray*}

Similarly, we can obtain
\begin{eqnarray*}
Z_jf(z)&=&\frac{{1}}{2}\left(\frac{{1}}{\rho}\frac{\partial{\tilde{a}}}{\partial{\rho}}-\frac{{1}}{2}\tilde{a}\right)
P_0\\&&+\left[\left\{\frac{{1}}{2(n+p+q-1)}\left(\rho\frac{\partial}{\partial{\rho}}-\frac{{1}}{2}\rho^2
\right)+1\right\}\tilde{a}\right]\frac{\partial{P}}{\partial z_j}.
\end{eqnarray*}
Hence the projection $\Pi_{p,q-1}$ of $\bar{Z_j}f$ is
\begin{equation}
\Pi_{p,q-1}(\bar{Z_j}f)=\left[\left\{\frac{{1}}{2(n+p+q-1)}\left(\rho\frac{\partial}{\partial{\rho}}+\frac{{1}}{2}\rho^2
\right)+1\right\}\tilde{a}\right]\frac{\partial{P}}{\partial{\bar{z_j}}}.
\end{equation}
Let $p=0$ and $q=1.$ Then there exists a $j$ such that
$\dfrac{\partial{P}}{\partial{\bar{z_j}}}$ is a non-zero constant.
Therefore, in this case,

\[\Pi_{0,0}(\bar{Z_j}f)(z)=C\left\{\frac{1}{2n}\left(\rho\frac{\partial}{\partial{\rho}}+
\frac{1}{2}\rho^2\right)+1\right\}\tilde{a}(\rho)\] is in
$Z(Ann(r, R)).$ Evaluating the twisted spherical mean at $z=0,$ we
get

\[\left\{\frac{1}{2n}\left(\rho\frac{\partial}{\partial{\rho}}+
\frac{1}{2}\rho^2 \right)+1\right\}\tilde{a}=0 .\]

To solve this equation, we substitute
$\tilde{a}(\rho)={e^{-\frac{1}{4}\rho^2}}\tilde{b}(\rho)$ and get
the differential equation
\[e^{-\frac{1}{4}\rho^2}\left\{\frac{1}{2n}\rho\frac{\partial}{\partial{\rho}}+1\right\}\tilde{b}=0.\]

Solving it, we conclude that for the case $p=0,q=1$, the coefficient
$\tilde{a}(\rho)=c_1{e^{-\frac{1}{4}\rho^2}}\rho^{-2n}.$

A simple induction argument gives that for $p=0$ and $q \geq 1,$
$\tilde{a}$ satisfies

\[\prod_{i=1}^q\left\{\frac{{1}}{2(n+q-i)}\left(\rho\frac{\partial}{\partial{\rho}}+\frac{{1}}{2}\rho^2
\right)+1\right\}\tilde{a}=0\] and therefore
\[~\tilde{a}(\rho)=~\sum_{i=1}^q~c_{i}~e^{-\frac{1}{4}\rho^2}\rho^{-2(n+q-i)}.\]

Similarly, using equation (3.12), we find that for $p \geq 1,$ and
$q=0,$ we have
\[~\tilde{a}(\rho)=~\sum_{k=1}^p~d_{k}~e^{\frac{1}{4}\rho^2}\rho^{-2(n+p-k)}.\]
This completes the description of the coefficients $(p,q)$ when
either $p$ or $q$ is zero.

Next we take up the case when $p=1, q=1.$ This can be reduced to
case of $p=0, q=0,$ by means of the operators $Z_j $ and
$\bar{Z_j}.$

For this, using Lemma \ref{lemma1}, without loss of generality,
assume that function is of the form, $f(z)=\tilde{a}(\rho)z_1
\bar{z_2}\in Z(Ann(r, R)).$ Applying the operators $Z_1\bar{Z_2}$
and taking the $(0,0)^{th}$ projection, we have
\[\left\{\frac{{1}}{2(n+1)}\left(\rho\frac{\partial}{\partial{\rho}}-\frac{{1}}{2}\rho^2
\right)+1\right\}\left\{\frac{{1}}{2(n+1)}\left(\rho\frac{\partial}{\partial{\rho}}+
\frac{{1}}{2}\rho^2\right)+1\right\}\tilde{a}=0. \] Solving this
differential equation, we get
\[\tilde{a}(\rho)=c_1{e^{\frac{1}{4}\rho^2}}\rho^{-2(n+1)}+d_1{e^{-\frac{1}{4}\rho^2}}\rho^{-2(n+1)}.\]

Finally, for the arbitrary $p,q$, again using Lemma \ref{lemma1}, we can
again assume that the function is of the form,
$f(z)=\tilde{a}(\rho)z_1^p\bar{z_2}^q \in Z(Ann(r,R)).$

Applying the operator $Z_1^p\bar{Z_2^q}$ and taking $(0,0)^{th}$
projection, we have
\[\prod_{i=1}^p\left\{A_i\left(\rho\frac{\partial}{\partial{\rho}}-\frac{1}{2}\rho^2
\right)+1\right\}\prod_{k=1}^q\left\{B_k\left(\rho\frac{\partial}{\partial{\rho}}+
\frac{1}{2}\rho^2\right)+1\right\}\tilde{a}=0, \] where $A_i=
{{(2(n+p+q-i))}}^{-1}$ and $B_k={(2(n+p+q-k))}^{-1}.$

 Solving
this, we get
\[~ \tilde{a}(\rho)=
~\sum_{i=1}^p~c_{i}~e^{\frac{1}{4}\rho^2}\rho^{-2(n+p+q-i)}+
~\sum_{k=1}^q~d_{k}~e^{-\frac{1}{4}\rho^2}\rho^{-2(n+p+q-k))}\]

This completes the proof of the theorem.

\end{proof}

Now we shall prove the sufficient part of Theorem ~\ref{th1}. The
proof of this part runs exactly the same way as that worked out
for an example in \cite{NT1}. Nonetheless, for the sake of completeness,
we give it here for the general case. This proof will be using the
result of Epstein and Kliener \cite{EK} on the spherical means on
$\mathbb R^d,$ which we briefly describe here.

For a function $f$ on $\mathbb R^d$ we have the spherical harmonic
expansion $$ f(x) = f( \rho \omega) = \sum_{k=0}^{\infty}
\sum_{l=1}^{d_k}~a_{kl}(\rho)~Y_k^l(\omega)$$ where $\rho = |x|$
and $\{ Y_k^l(\omega):~l = 1, 2, \cdots \cdots d_k \}$ is an
orthonormal basis for the space $V_k$ of homogeneous harmonic
polynomials of degree $k$ restricted to the unit sphere. For each
$k,$ the space $V_k$ is invariant under the action of $SO(d).$
When $d = 2m$ for some $m,$ it is invariant under the the action
of the unitary group $U(m)$ as well, and under this action of
$U(m)$ the space $V_k$ breaks up into an orthogonal direct sum of
$H_{p,q}$'s where $p+q = k.$
 Let $\sigma_s$ stand for the normalized surface measure on the
sphere of radius $s$ centered at the origin contained in $\mathbb
R^d.$ The main result in \cite{EK} implies the following theorem for the special case of the annulus
$\{x\in \mathbb R^d :|x|>B\}$:

\bt\label{th5} A continuous function $f$ on $\R ^d$ satisfies
$$
\int_{|y| =s }g(x+y)d\sigma_s(y) = 0  \quad \mbox{for} \quad s
> |x| + B \quad \mbox{for all} \quad x  \in \R^d $$ if and only if
$$a_{kl}(\rho) = \sum_{i=0}^{k-1}~\alpha_{kl}^{i}~\rho^{k-d-2i},~\quad \alpha_{kl}^{i} \in ~\mathbb C
,$$ for all $k > 0,$ $1 \leq l \leq d_k,$ and $a_0(\rho) = 0$
whenever $\rho > B.$ \et

Next we take up the proof of the sufficient part of Theorem
\ref{th1}.
 \bt\label{th6}
Suppose $h$ is a function defined on $Ann(r,\infty)$ by $
h(z)=\dfrac{{e^{\frac{1}{4}|z|^2}}P(z)}{|z|^{2(n+p+q-i)}},$ where
$P\in H_{p,q}$ and $1\leq i\leq p.$ Then $h \in Z(Ann(r,
\infty)).$ \et

\begin{proof} We have to show that
$ h\times \mu_s (z) = 0 $ for all $z, s $ with $|z|+r<s.$

Consider, $$ h \times \mu_s(z) = \int_{|w| =s}
~\frac{{e^{\frac{1}{4}|z+w|^2}}P(z+w)}{|z+w|^{2(n+p+q-i)}}e^{-\frac{i}{2}Im(z.\bar{w})}
d\mu_s(w).$$

Expanding the term $|z+w|^2$ and simplifying, we see that it is
enough to consider the integral
$$ \int_{|w| = s
}~\frac{{e^{\bar{z}.w}}P(z+w)}{|z+w|^{2(n+p+q-i)}} d\mu_s(w).$$

On expanding the exponential factor,  this leads to terms of the
form
$$\int_{|w| = s}~\frac{
{w^\alpha}P(z+w)}{|z+w|^{2(n+p+q-i)}}d\mu_s(w)$$ where $\alpha$ is
a multi-index. Writing ${w_1} = z_1+w_1 - z_1$ etc. and expanding
again we see that it is enough to consider terms of the form

 \[
\int_{|w| =s }~\frac{(w+z)^\beta P(z+w)}{|z+w|^{2(n+p+q-i)}}
d\mu_s(w).\] Let $g(z)=\frac {z^\beta P(z)}{|z|^{2(n+p+q-i)}}.$
Then, the above expression is

 $$
\int_{|w| =s }g(z+w)d \sigma_s(w),$$ which is a Euclidean
spherical mean of $g$ on the sphere of radius $s$ centered at the
origin contained in $\mathbb R^{2n}.$

Thus, we need to show that $$ \int_{|w| =s }g(z+w)d
\sigma_s(w)=0,$$ for $s>|z|+r.$

Using the Lemma \ref{lemma3}, we have the decomposition
$$z^\beta P(z)= P_0(z)+|z|^{2}P_1(z)+ \dots+|z|^{2l}P_l(z)$$ where $P_j \in
H_{p+|\beta|-j,q-j}$, ~for $0\leq j \leq l$ , ~$l\leq
min(|\beta|,q).$
\medskip
With this, the function $g$ further decomposes in functions of the
form $|z|^{- 2(n+p+q-i-j)}P_j(z)$.

Hence, to prove that $g$ satisfies the desired convolution
equation, it is enough to show that the function $|z|^{-
2(n+p+q-i-j)}P_j(z)$ satisfies it.

Let us rewrite $$|z|^{-
2(n+p+q-i-j)}P_j(z)=\rho^{k-2n-2(p+q-i-j)}Y_k,$$ where
$k=p+q+|\beta| -2j$ and $Y_k$ is a spherical harmonic of degree
$k$ on $\mathbb{R}^{2n}$. Using Theorem \ref{th5}, we need to show
that $ 0\leq p+q-i-j\leq k-1$, or equivalently $j-i\leq|\beta|-1.$

If $|\beta|\leq q$, ~then $l=|\beta|$ and $ j-i\leq
j-1\leq|\beta|-1$ (since $j\leq |\beta|$ and $1\leq i\leq p$). For
$|\beta|> q$, we get $l=q$. Since $|\beta|> q\geq j$, therefore we
have $|\beta|-1>j-1\geq j-i$, as $1\leq i\leq p$.

This completes the proof.
\end{proof}

Similarly, we can prove the following theorem.

\bt \label{th7} Suppose $h$ is a function defined on $Ann(r,R)$ by
$ h(z)=\dfrac{{e^{-\frac{1}{4}|z|^2}}P(z)}{|z|^{2(n+p+q-k)}},$
where $P\in H_{p,q}$ and $k=1,\cdots,q.$ Then $h \in Z(Ann(r,
R)).$ \et

Putting together Theorem \ref{th6} and Theorem \ref{th7}, the
sufficient part of Theorem \ref{th1} follows.
\medskip

\section{\bf Proofs of the support theorems and concluding remarks}

We begin by recalling the Helgason's support theorem (\cite{H}, p. 16) for Euclidean spherical means.

\bt\label{th81} Let $g$ be a continuous function on $\R ^d$ such that for each $k=0,1,\cdots,$
${\mbox sup}\,|x|^k|g(x)|<\infty.$ Then $g$ is supported in $\{x \in \mathbb R^d: |x| \leq B\}$ if and only if
$$
\int_{|y| =s }g(x+y)d\sigma_s(y) = 0  \quad \mbox{for} \quad s
> |x| + B \quad \mbox{for all} \quad x  \in \R^d. $$
\et

Here, as before, $\sigma_s$ stand for the normalized
surface measure on $\{x \in \mathbb R^d : |x| = s \}.$

This theorem can now be deduced as a corollary of Theorem \ref{th5} (also noted in  \cite{EK}), as the spherical harmonic coefficients of $f$
satisfy the same decay conditions as $f$.

Next we recall
the following support theorems for the twisted spherical
means proved in [NT2] for the twisted spherical
means.

\bt Let $f$ be a function on $\mathbb C^n$ such that $f(z)
e^{\frac{1}{4}|z|^2}$ is in the Schwartz class. Then $f$ is
supported in $|z| \leq r$ if and only if $f \times \mu_s(z) = 0 $
for $s > r + |z|$ for every $z \in \mathbb C^n.$ \et

In the above theorem, the function $f$ is assumed to have exponential decay, which reflects the non-Euclidean nature of the
twisted spherical means. Such decay conditions also arise naturally in the integral geometry on the Heisenberg group as can be seen in the results in \cite{AR}, \cite{NT1}. However, the differentiability conditions on the function are
genuine and cannot be relaxed.
This is because the condition that $f(z) e^{\frac{1}{4}|z|^2}$ is in Schwartz class is not translation invariant (\cite{NT2}). Nonetheless, to do away with the smoothness condition on $f,$ a stronger condition like $|f(z)| \leq C~e^{-(\frac{1}{4}+\epsilon)|z|^2},$ for some $\epsilon >0$ can be imposed. As then we may convolve $f$ on the right with a radial approximate
identity to get smooth functions $\{f_\epsilon\}$ which approximate $f$ and also satisfy the vanishing mean conditions.

In contrast, in Theorem \ref{th2} we do not impose any differentiability conditions on the function nor do we impose a stronger decay condition. Our conditions can be thought of as an exact analogue of the conditions in the Euclidean set up.

The proof of Theorem \ref{th2} follows immediately from Theorem \ref{th1}, as the spherical harmonic coefficients $a_j^{p,q}$ satisfy the same decay conditions as the function $f$.

When $n = 1,$ the authors in \cite{NT2} have shown that under very weak conditions on $f$ and with a suitable condition involving both sided twisted spherical
means the following result holds.

\bt Let $f$ be a locally integrable function on $\mathbb C$ satisfying the condition that $ |f(z)| \leq ~C~e^{\frac{1}{4}(1-\epsilon)|z|^2}
$ for some $\epsilon > 0.$ Then $f$ is supported in $|z| \leq B$
if and only if $f \times \mu_r (z) = \mu_r \times f(z)  = 0 $ for
$ r > B+ |z|$ for every $z \in ~\C.$\et

In the version Theorem \ref{th3} of this support theorem, we do not need any growth conditions on the function.

For a proof of Theorem \ref{th3}, let us consider the space $Z^\ast(Ann(r, R))$ of continuous
functions $f$ on $~\mathbb C^n$ with both the twisted spherical
means $f \times \mu_s (z) = \mu_s \times f(z) = 0 $~~for all
spheres $S_s(z)$ contained in the $Ann(r, R) $ and with
$B_r(0)\subseteq B_s(z)$. Then Theorem \ref{th1} can be strengthened to the following result:

 \bt\label{th8} A necessary and sufficient condition for a function $f$ to belong to $Z^\ast(Ann(r,
R))$ is that for $p, q \in \mathbb Z^+$ and $1\leq j\leq d(p,q),$
the spherical harmonic coefficients $~~\tilde{a}_j^{p,q}(\rho)$ of $f$ satisfy, for $r<\rho<R,$

$$~\tilde{a}_j^{p,q}(\rho)=
\sum_{i=1}^{min(p,q)}~c_i~e^{\frac{1}{4}\rho^2}\rho^{-2(p+q+n-i)}+
\sum_{k=1}^{min(p,q)}~d_k~e^{-\frac{1}{4}\rho^2}\rho^{-2(p+q+n-k)},
p \neq 0, q \neq 0$$ and $~\tilde{a}_j^{0,q}=\tilde{a}_j^{p,0}=0.$
Here $c_i, d_k$ are arbitrary constants in $\mathbb C.$\et

\begin{proof} As
$ \mu_\rho \times f = \overline{\bar{f} \times \mu_\rho},$ it follows that $f\in Z^\ast(Ann(r, R))$ if and only if $\bar{f}\in
Z^\ast(Ann(r, R)).$ Also a $(p,q)$th spherical harmonic coefficient
of $f,$ $\tilde{a}_j^{p,q}(f)$ is related to the corresponding
spherical harmonic coefficient of $\bar{f}$ by
$\tilde{a}_j^{p,q}(\bar{f})=\overline{\tilde{a}_j^{p,q}(f)}.$ Hence
the conclusion follows from Theorem ~\ref{th1}.
\end{proof}

The proof of Theorem \ref{th3} now follows as a corollary of the above theorem and the observation that for $n=1$, the nonzero spaces $H_{p,q}$ will have either
the $p=0$ or $q=0.$

It is therefore no surprise that the
decay condition on $f$ could be completely relaxed for the support theorem on functions on $\mathbb C.$

\vskip.1in

Finally, coming back to the Heisenberg group $\mathbb H^n=\mathbb C^n\times \mathbb R,$ let
$f$ be a continuous function on $\mathbb H^n$ which has the spherical means (as defined in \ref{exp2})
$f\ast\mu_s(z,t)=0$ for all $t\in \mathbb R$ and $z\in \mathbb C^n$ satisfying
$B_r(0)\subseteq B_s(z)$ and
 $S_s(z)\subseteq Ann(r,R).$
 The problem of characterizing such functions in general is open.
 However, if $f$ is of the form
 $f(z,t)=e^{i\lambda t}\varphi(z),~\lambda\in \mathbb R\setminus\{0\},$
 then an easy modification of the proof of Theorem \ref{th1} for $\lambda$-twisted spherical means,
 $\lambda$ in $\mathbb R\setminus\{0\},$
gives a  characterization for $f$ in terms of the spherical harmonic
 coefficients of the function $\varphi$. For $\lambda=0$, the problem reduces to the problem on Euclidean spherical means.

\vspace{.1in}

\noindent {\bf Acknowledgements:}~The authors wish to thank Dr. E.
K. Narayanan for several fruitful discussions. We also thank the referee for some useful comments. The second author
wishes to thank the MHRD, India, for the senior research
fellowship and IIT Kanpur for the support provided during the
period of this work.

\vskip.15in
\begin{flushleft}
Department of Mathematics and Statistics, \\
Indian Institute of Technology \\
Kanpur 208 016 \\
India\\
E-mail:~rrawat@iitk.ac.in, rksri@iitk.ac.in\\

\end{flushleft}

\end{document}